\documentclass[12pt]{article}
\usepackage[left=3cm,top=3.0cm,right=3cm,bottom=2.6cm]{geometry}

\usepackage[ansinew]{inputenc}
\usepackage{amscd,amsfonts,amsmath,amssymb,amsthm,enumerate,epsfig,float,graphics,graphicx,pst-all,yhmath}
\newtheorem{lemma}{Lemma}

\newtheorem{theorem}{Theorem}
\newtheorem{corollary}{Corollary}

\newtheorem{remark}{Remark}

\newcommand{\fin}{\hfill $\Box$}

\title {Some results about equichordal convex bodies}
\author{
Jes\'us Jer\'onimo-Castro$^{1}$,  Francisco G. Jimenez-Lopez$^{2}$,\\ and Efr\'en Morales-Amaya$^{3}$\\ 
\small{$^{1,2}$Facultad de Ingenier\'ia}\\
\small{Universidad Aut\'onoma de Quer\'etaro, M\'exico}\\
\small{$^{3}$Facultad de Matem\'aticas-Acapulco,}\\
\small{Universidad Aut\'onoma de Guerrero, M\'exico}\\
\small{$^{1}$jesus.jeronimo@uaq.mx,}\\
\small{$^{2}$fjimenez@uaq.mx,}\\
\small{$^{3}$emoralesamaya@gmail.com}
}
\begin{document}

\maketitle

\begin{abstract}
Let $K$ and $L$ be two convex bodies in $\mathbb R^n$, $n\geq 2$, with $L\subset \text{int}\, K$. We say that $L$ is an \emph{equichordal body} for $K$ if every chord of $K$ tangent to $L$ has length equal to a given fixed value $\lambda$. In \cite{Barker}, J. Barker and D. Larman proved that if  $L$ is a ball, then $K$ is a ball concentric with $L$. In this paper we prove, derived from the proof of Theorem \ref{cuerdas}, that there exist an infinite number of closed curves, different from circles, which possess an equichordal convex body. If the dimension of the space is more than or equal to 3, then only Euclidean balls possess an equichordal convex body. We also prove some results about isoptic curves and give relations between isoptic curves and convex rotors in the plane.
\end{abstract}

\section{Introduction}
Let $K$ be a convex body in the plane, i.e., a compact and convex set with non-empty interior, and let $\mathcal P$ be a convex polygon. It is said that $K$ is a rotor in $\mathcal P$ if for every rotation $\rho$, there is a translate of $\mathcal P$ that contains to $\rho(K)$ and all sides of $\mathcal P$ are tangent to $K$. There are many results about rotors in regular polygons, see for instance \cite{Groemer}, and for the particular case of rotors in equilateral triangles see \cite{Yaglom}. The case of rotors in squares is well-known, indeed, bodies of constant width are a very important topic in Convex Geometry and have many interesting properties and applications in mechanisms (see the quite nice book \cite{Martini}). 
 
Another topic, apparently not related to rotors, is the \emph{Equichordal Problem}. Let $x$ be a point in the interior of a convex body $K$, we say that $x$ is an \emph{equichordal point} if every chord of $K$ through $x$ have the same length. The famous Equichordal Problem, due to W. Blaschke, W. Rothe, and R. Weitzenb\"ock \cite{Blaschke},  asks about the existence of a convex body with two equichordal points. There are many false proofs about the non existence of such a body, however,  M. Rychlik finally gave a complete proof about the non existence of a body with two equichordal points in \cite{Rychlik}. It is worth mentioning that there are many convex bodies, different from the disc, which have exactly one equichordal point. Here we are interested in a generalization of the notion of equichordal point in the following way: Let $K$ and $L$ be two convex bodies in $\mathbb R^n$, $n\geq 2$, with $L\subset \text{int}\, K$. We say that $L$ is an \emph{equichordal body} for $K$ if every chord of $K$ tangent to $L$ have length equal to a given fixed value $\lambda$. In \cite{Barker}, J. Barker and D. Larman proved that if $K$ is a convex body that possesses an equichordal ball then it is also a ball. However, we wonder if there exist convex bodies different from balls which possess an equichordal convex body in its interior. It seems that bodies which float in equilibrium in every position provide examples of such bodies in the plane (see for instance \cite{Wegner}), however, it is not clear if the considered bodies $K$ are convex or not.  

In sections 3 and 4 of this paper we study convex bodies $L$ for which the chords of one of its isoptic curves (defined in the following section), that are tangent to $L$ have length equal to a constant number $\lambda$. Moreover, these bodies $L$ are examples of rotors in regular polygons and if we fix the convex body $L$ and the circumscribed polygon is rotated, while maintained circumscribed to $K$, its vertices describe the isoptic curve of $L$. In section 6 we also prove that in dimension $3$ or higher, only Euclidean balls (or simply balls) possess and equichordal convex body in its interior.

%%%%%%%%%%%%%%%%%%%%%%%%%%%%%%%%%%%%%%%%%%%%%%%%%%%%%%%%%%%%%%%%%%%%%%%%%%%%%%%%%%%%O

\section{Preliminary concepts}
We give first some definitions and notation. Let $K$ be a given planar convex body; for every real number $t$ we denote by $\ell(t)$ the support line of $K$ with outward normal vector $u(t)=(\cos t, \sin t)$. The function $p:\mathbb R\longrightarrow \mathbb R$, defined as $p(t)=\max_{x\in K} \langle u(t), x\rangle$, is known by the name of support function of $K$. When the origin $O$ is contained in $K$, $p(t)$ is nothing else than the distance from $O$ to the support line $\ell(t).$ The distance between the support lines $\ell(t)$ and $\ell(t+\pi)$  is called the width of $K$ in direction $u(t)$ and it is denoted by $w(t)$, in other words, $w(t)=p(t)+p(t+\pi)$. If $w(t)$ is constant, independently of $t$, we say that $K$ is a body of constant width. For any $\alpha\in (0,\pi)$, the $\alpha$-isoptic $K_{\alpha}$ of $K$ is
defined as the locus of points at which two tangent lines to $K$
intersect at an angle $\alpha$.  Using the support function, $\partial K$ is parameterized (see for instance \cite{Val})
by 
\begin{equation}{\label{eq1}}
\gamma(t)=p(t)u(t)+ p'(t) u'(t), \ \text{for} \  t\in[0,2\pi].
\end{equation}

The isoptic curve $K_{\alpha}$ can be parameterized by
the same angle by the formula (see \cite{Cies_Mier_Moz} or \cite{Jeronimo})
\begin{equation}{\label{eq2}}
\gamma_{\alpha}(t)=p(t)u(t)+\left[p(t)\cot\alpha+\frac{1}{\sin\alpha}p(t+\pi-\alpha)\right]u'(t).
\end{equation}

By Cauchy's formula, the perimeter of $K$ can be obtained by (see \cite{Val})
\begin{equation}\label{Cauchy}
L(K)=\int_0^{2\pi}p(t)dt.
\end{equation}

\begin{figure}[H]
    \centering
    \includegraphics[width=1.0\textwidth]{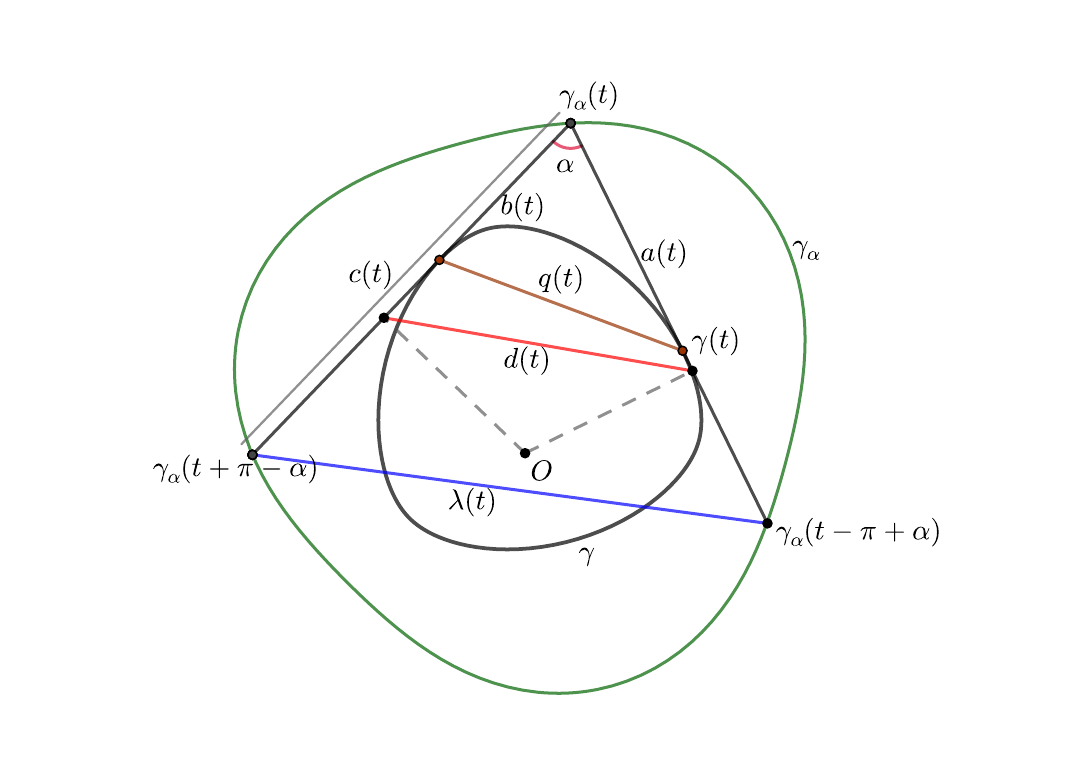}
    \caption{Parameters of the isoptic curve}
    \label{parametros}
\end{figure}

For any $t\in\mathbb R$ we define (see Fig. \ref{parametros})
\begin{align}
    a(t)&=|\gamma_{\alpha}(t)-\gamma(t)|,\\
    b(t)&=|\gamma_{\alpha}(t)-\gamma(t+\pi-\alpha)|,\\
    c(t)&=|\gamma_{\alpha}(t)-\gamma_{\alpha}(t+\pi-\alpha)|=b(t)+a(t+\pi-\alpha),\\
    q(t)&=|\gamma(t)-\gamma(t+\pi-\alpha)|,\\
    \lambda(t)&=|\gamma_{\alpha}(t+\pi-\alpha)-\gamma_{\alpha}(t-\pi+\alpha)|.
\end{align}

We also define $d(t)$ to be the distance between the points obtained by projecting the origin $O$ onto the support lines of $K$ at $\gamma(t)$ and $\gamma(t+\pi-\alpha)$.

By some tedious but simple calculations we can express the lengths of all these chords in terms of the support function of $K$: 
\begin{align}
        a(t)&=\frac{1}{\sin\alpha}[p(t+\pi-\alpha)+p(t)\cos\alpha-p'(t)\sin\alpha]\label{a_t},\\
        b(t)&=\frac{1}{\sin\alpha}[p(t+\pi-\alpha)\cos\alpha+p'(t+\pi-\alpha)\sin\alpha+p(t)],\label{b_t}\\
        c(t)&=\frac{1}{\sin\alpha}[2p(t+\pi-\alpha)\cos\alpha+p(t)+p(t-2\alpha)],\label{c_t}\\
        d(t)&=\sqrt{p^2(t)+p^2(t+\pi-\alpha)+2p(t)p(t+\pi-\alpha)\cos\alpha}.\label{d_t}  
    \end{align}

%\begin{figure}[H]
%    \centering
%    \includegraphics[width=.8\textwidth]{tangente_isoptica.pdf}
%    \caption{ $a(t)=\frac{1}{\sin\alpha}[p(t+\pi-\alpha)+p(t)\cos\alpha-p'(t)\sin\alpha]$}
%    \label{tangente_isoptica}
%\end{figure}

Finally, the support function of a convex body is a periodic function, with period $2\pi$, and it is also absolutely continuous, so we can consider its expansion in terms of the Fourier series (see \cite{Groemer}), i.e.,
    \begin{equation}\label{fourier1}
        p(t)=a_0+\sum_{n=1}^\infty\left(a_n\cos nt+b_n\sin nt\right).
    \end{equation}
    
    The first and second derivatives of $p$ are expressed as
\begin{align}
        p'(t)&=-\sum_{n=1}^\infty\left(na_n\sin nt-nb_n\cos nt\right),\label{fourier2}\\
        p''(t)&=-n^2 \sum _{n=1}^{\infty }\left(a_n \cos n t+b_n \sin n t\right).\label{fourier3}
\end{align}

%%%%%%%%%%%%%%%%%%%%%%%%%%%%%%%%%%%%%%%%%%%%%%%%%%%%%%%%%%%%%%%%%%%%%%%%%%%%%%%%%%

%%%%%%%%%%%%%%%%%%%%%%%%%%%%%%%%%%%%%%%%%%%%%%%%%%%%%%%%%%%%%%%%%%%%%%%%%%%%%%%%%%
\section{Some results about isoptic curves in the plane}
Our first result about isoptic chords is the following.

\begin{theorem}\label{cuerdas}
Let $K$ be a strictly convex body in the plane with differentiable boundary and let $\alpha\in(0,\pi)$ be a fixed angle such that $\frac{\alpha}{\pi}$ is an irrational number. Suppose $c(t)=c_0$, for every $t\in [0,2\pi]$, for a positive number $c_0.$ Then $K$ is a disc.
\end{theorem}

\emph{Proof.} Since $c(t)=c_0$ we have by (\ref{c_t}) that $$c(t)=\frac{1}{\sin\alpha}[2p(t+\pi-\alpha)\cos\alpha+p(t)+p(t-2\alpha)]=c_0,$$ it follows that \begin{equation}\label{dif1}
2p'(t+\pi-\alpha)\cos\alpha+p'(t)+p'(t-2\alpha)=0.
\end{equation}
If we substitute the Fourier coefficients of $p(t)$ in the differential equation (\ref{dif1}), by (\ref{fourier2}) we have
\begin{align*}
2\cos \alpha & \sum_{n=1}^{\infty}(nb_n\cos n(t+\pi-\alpha)-na_n\sin
n(t+\pi-\alpha))\\
+ & \sum_{n=1}^{\infty}(nb_n\cos nt-na_n\sin
nt)\\
+ & \sum_{n=1}^{\infty}(nb_n\cos n(t-2\alpha)-na_n\sin
n(t-2\alpha))=0.
\end{align*}  
Since this holds for every real number $t$, we must have that for every $n$ 
\begin{align*}
& \cos nt[-2na_n \sin n(\pi-\alpha)\cos\alpha +n a_n\sin 2n\alpha +2nb_n \cos n(\pi-\alpha)\cos \alpha \\
& +n b_n +n b_n\cos 2n\alpha] + \sin nt[-2na_n \cos n(\pi-\alpha)\cos\alpha - n a_n -na_n \cos 2n\alpha \\
&-2nb_n \sin n(\pi-\alpha)\cos\alpha +n b_n\sin 2n\alpha]=0.
\end{align*}
The coefficients of $\cos nt$ and $\sin nt$ must be both equal to 0, hence we have that
$$\begin{bmatrix}
-2\sin n(\pi-\alpha)\cos\alpha +\sin 2n\alpha & 2\cos n(\pi-\alpha)\cos\alpha+1 +\cos 2n\alpha\\
-2\cos n(\pi-\alpha)\cos\alpha-1-\cos 2n\alpha & -2\sin n(\pi-\alpha)\cos\alpha +\sin 2n\alpha
\end{bmatrix}\begin{bmatrix}a_n\\b_n \end{bmatrix}=\begin{bmatrix}0\\ 0 \end{bmatrix}.$$
The determinant of the matrix above is given by $$(-2\sin n(\pi-\alpha)\cos\alpha +\sin 2n\alpha)^2+(2\cos n(\pi-\alpha)\cos\alpha+1 +\cos 2n\alpha)^2.$$
This determinant is zero only if 
\begin{equation}\label{launo}
-2\sin n(\pi-\alpha)\cos\alpha +\sin 2n\alpha=0
\end{equation}
and 
\begin{equation}\label{lados}
2\cos n(\pi-\alpha)\cos\alpha+1 +\cos 2n\alpha=0.
\end{equation}
Since $\alpha\in (0,\pi)$, for every $n\geq 2$ we have that none of (\ref{launo}) and (\ref{lados}) are satisfied if $\frac{\alpha}{\pi}$ is an irrational number. Hence, we have that the determinant is non zero and then $a_n=b_n=0$ for every $n\geq2$. It follows that $p(t)=a_0+a_1\cos t+b_1 \sin t$, i.e., $p$ is the support function of a disc (see for instance \cite{Val}).
\fin

\begin{figure}[H]
    \centering
    \includegraphics[width=.66\textwidth]{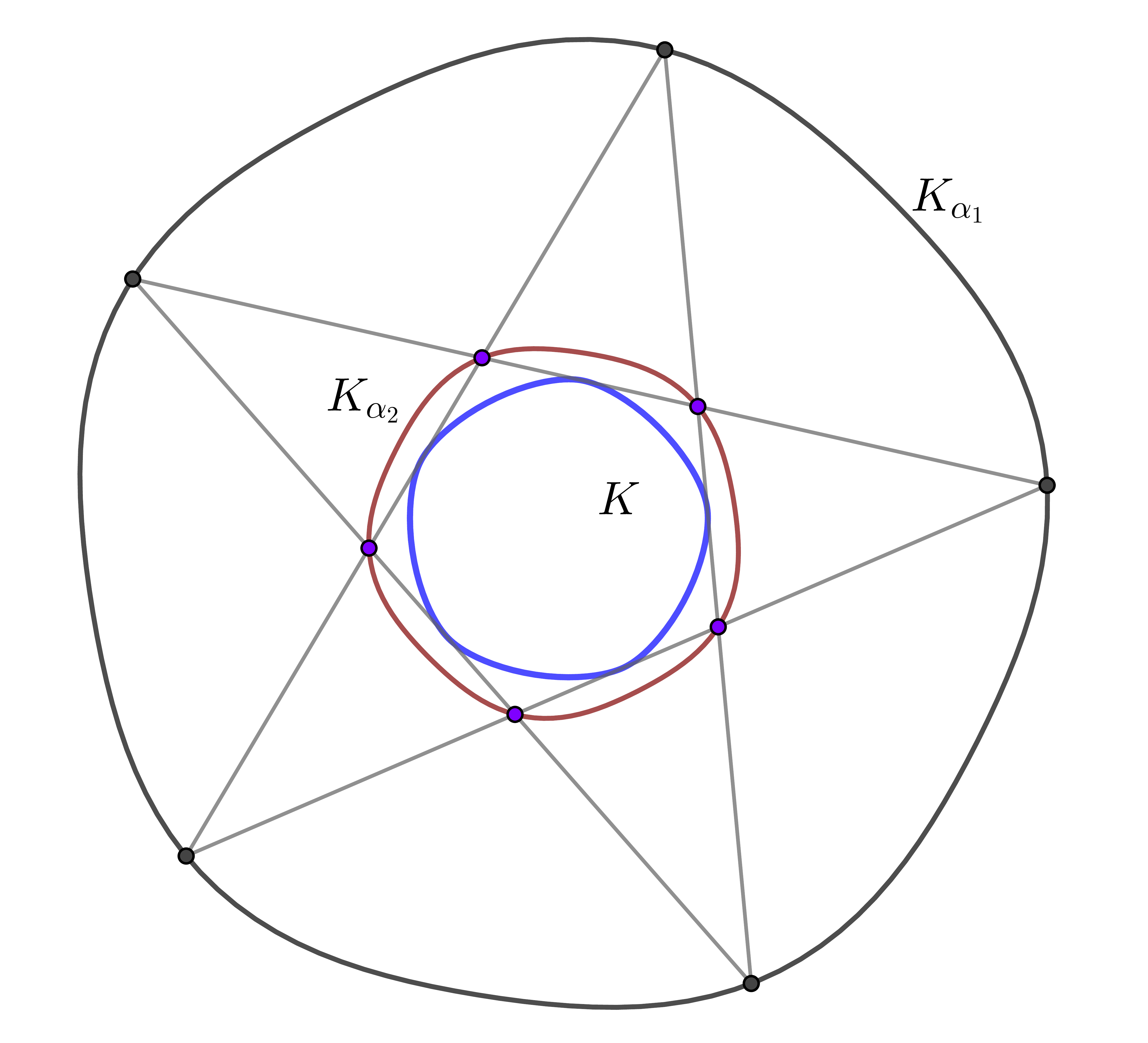}
    \caption{A convex body with two isoptics with corresponding values of $c(t)$ both of constant value for the angles $\alpha_1=\frac{\pi}{5}$ and $\alpha_2=\frac{3\pi}{5}$}
    \label{rotor_pentagono}
\end{figure}

\begin{remark}
If $\frac{\alpha}{\pi}$ is a rational number, then there exist convex bodies $K$ different from discs for which $c(t)$ is constant. For instance: for the angles $\alpha_1=\frac{\pi}{5}$ and $\alpha_2=\frac{3\pi}{5}$ let $K$ be the convex body whose support function is given by $p(t)=60+\cos 5t+\sin 5t$ (see Fig. \ref{rotor_pentagono}). In this case the isoptic curves $K_{\alpha_1}$ and $K_{\alpha_2}$ have the property that its chords tangent to $K$ have constant values $c_1$ and $c_2$, respectively, with $\frac{c_1}{c_2}=\frac{1}{2-\tau}$, where $\tau=\frac{1+\sqrt 5}{2}$. Moreover, $K_{\alpha_1}$ is homothetic to $K_{\alpha_2}$ with ratio of homothety equal to $-\frac{1}{2-\tau}$.
\end{remark}

However, if we impose the additional condition that $q(t)$ is also constant, then $K$ must be a disc.

\begin{theorem}\label{cuerdas_mas_q}
Let $K$ be a strictly convex body in the plane with differentiable boundary $\partial K$, and let $\alpha\in(0,\pi)$ be a fixed angle. Suppose $c(t)=c_0$, and $q(t)=q_0$, for every $t\in [0,2\pi]$, and for two positive numbers $c_0$ and $q_0$. Then $K$ is a disc. 
\end{theorem}

\begin{figure}[H]
    \centering
    \includegraphics[width=.63\textwidth]{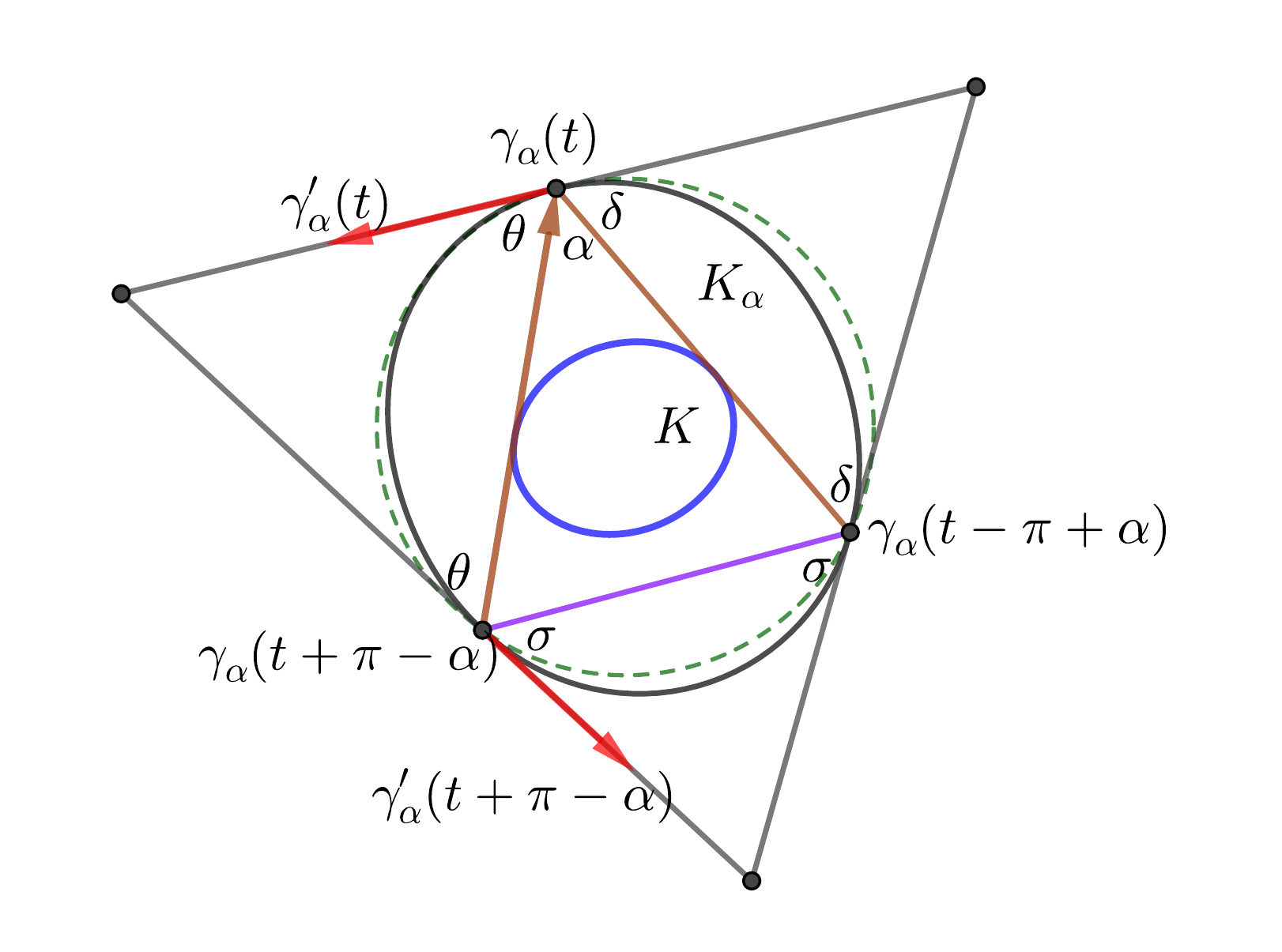}
    \caption{The angle between the chord $[\gamma_{\alpha}(t+\pi-\alpha),\gamma_{\alpha}(t)]$ and the tangent lines are equal}
    \label{vectores_tangentes}
\end{figure}

\emph{Proof.} By some simple calculations we have that $|\gamma'_{\alpha}(t)|=\frac{q(t)}{\sin\alpha}$ (see for instance \cite{Cies_Mier_Moz}).  Since $c(t)=c_0$ is also constant, we have $$\frac{d}{dt}(|\gamma_{\alpha}(t+\pi-\alpha)-\gamma_{\alpha}(t)|^2)=\frac{d}{dt}(\langle \gamma_{\alpha}(t+\pi-\alpha)-\gamma_{\alpha}(t),\gamma_{\alpha}(t+\pi-\alpha)-\gamma_{\alpha}(t) \rangle)=0,$$ hence $$\langle \gamma_{\alpha}(t+\pi-\alpha)-\gamma_{\alpha}(t),\gamma'_{\alpha}(t+\pi-\alpha)\rangle=\langle \gamma_{\alpha}(t+\pi-\alpha)-\gamma_{\alpha}(t),\gamma'_{\alpha}(t)\rangle.$$ It follows that $$c_0\cdot |\gamma'_{\alpha}(t+\pi-\alpha)|\cos\beta_1=c_0\cdot |\gamma'_{\alpha}(t)|\cos\beta_2,$$ which implies that $\beta_1=\beta_2$, where $\beta_1$ is the angle between the vectors $\gamma_{\alpha}(t+\pi-\alpha)-\gamma_{\alpha}(t)$ and $\gamma'_{\alpha}(t+\pi-\alpha)$, and $\beta_2$ is the angle between the vectors $\gamma_{\alpha}(t+\pi-\alpha)-\gamma_{\alpha}(t)$ and $\gamma'_{\alpha}(t)$. It follows that the angle between the chord $[\gamma_{\alpha}(t+\pi-\alpha),\gamma_{\alpha}(t)]$ and the tangent lines at $\gamma_{\alpha}(t)$ and $\gamma_{\alpha}(t+\pi-\alpha)$, are equal (see Fig. \ref{vectores_tangentes}). Similarly, we obtain that the angles between the chord $[\gamma_{\alpha}(t-\pi+\alpha),\gamma_{\alpha}(t)]$ and the tangent lines at $\gamma_{\alpha}(t)$ and $\gamma_{\alpha}(t-\pi+\alpha)$, are equal. Since the length of the tangent vector $\gamma'_{\alpha}(t)$ is constant for every $t$, and all the triangles $\triangle \gamma_{\alpha}(t-\pi+\alpha)\gamma_{\alpha}(t)\gamma_{\alpha}(t+\pi-\alpha)$ are congruent, we also have that the angles between the chord $[\gamma_{\alpha}(t-\pi+\alpha),\gamma_{\alpha}(t+\pi-\alpha)]$ and the tangent lines at $\gamma_{\alpha}(t+\pi-\alpha)$ and $\gamma_{\alpha}(t-\pi+\alpha)$, are equal. By elementary geometry we have that the circle circumscribed to $\triangle \gamma_{\alpha}(t-\pi+\alpha)\gamma_{\alpha}(t)\gamma_{\alpha}(t+\pi-\alpha)$ and the body $K_{\alpha}$, share the tangent lines at the points $\gamma_{\alpha}(t-\pi+\alpha)$, $\gamma_{\alpha}(t)$, and $\gamma_{\alpha}(t+\pi-\alpha)$; under this condition it was proved in Lemma 3.3 in \cite{Jeronimo_Jimenez} that $K_{\alpha}$ must be a disc. Now we use Theorem 2 (b) in \cite{Jeronimo_chords} and conclude that $K$ is a disc. \fin

Let $h(t)$ denote the length of the segment from
$\gamma_{\alpha}(t)$ to the projection of $O$ into the support line of $K$ at $\gamma(t)$. It is easy to show that
\begin{equation}\label{h(t)}
h(t)=\frac{1}{\sin\alpha}[p(t)\cos\alpha+p(t+\pi-\alpha)].
\end{equation}

\begin{figure}[H]
    \centering
    \includegraphics[width=.61\textwidth]{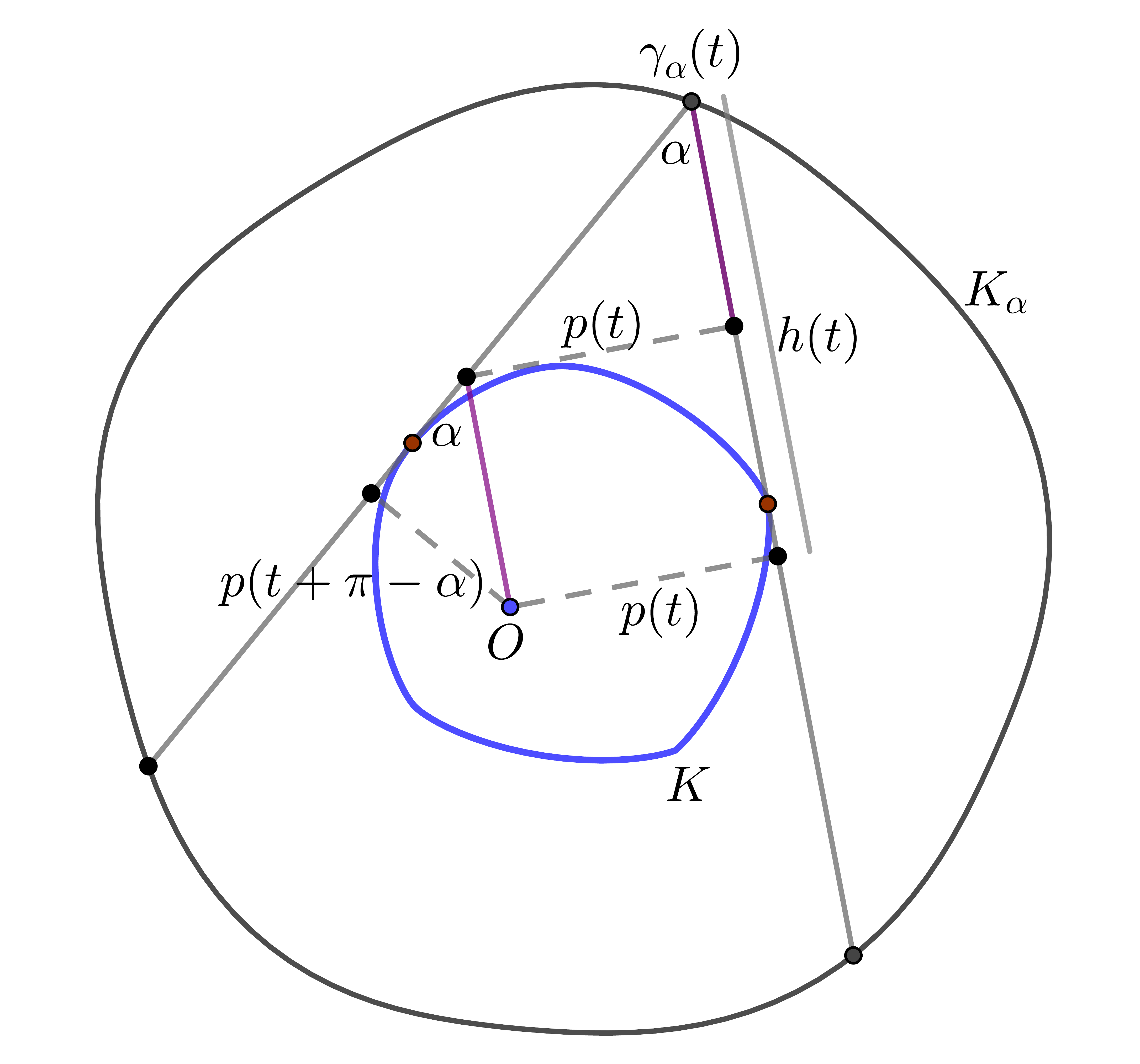}
    \caption{The value of $h(t)$ is constant}
    \label{parametro_h}
\end{figure}

We have the following result.

\begin{theorem}\label{parametroh}
Let $K$ be a strictly convex body in the plane with differentiable boundary and let $\alpha\in(0,\pi)$ be a fixed angle. Suppose $h(t)=h_0$, for every $t\in [0,2\pi]$, for a positive number $h_0.$ Then $K$ is a
disc.
\end{theorem}

\emph{Proof.} By (\ref{h(t)}) we have that
\begin{equation}\label{ecdif}
p'(t)\cos\alpha+p'(t+\pi-\alpha)=0.
\end{equation}

Using the expansion in Fourier series for $p'(t)$ (see equation (\ref{fourier2})) we have that
$$\sum_{n=1}^{\infty}(nb_n\cos nt\cos\alpha-na_n\sin
nt\cos\alpha+nb_n\cos n(t+\pi-\alpha)-na_n\sin n(t+\pi-\alpha))=0.$$
Using trigonometric identities and simplifying we conclude that
for all $n$ and for all $t\in[0,2\pi]$ we must have
\begin{align*}
0&=[-n\sin n(\pi-\alpha)a_n+(n\cos\alpha+n\cos
n(\pi-\alpha))b_n]\cos nt\\
&+[(-n\cos\alpha-n\cos n(\pi-\alpha))a_n-n\sin
n(\pi-\alpha)b_n]\sin nt.
\end{align*}
This yields to the system of equations
$$\begin{bmatrix}
-n\sin n(\pi-\alpha) & n\cos\alpha+n\cos n(\pi-\alpha)\\
-n\cos\alpha-n\cos n(\pi-\alpha) & -n\sin n(\pi-\alpha)
\end{bmatrix}\begin{bmatrix}a_n\\b_n \end{bmatrix}=\begin{bmatrix}0\\ 0 \end{bmatrix}.$$
Notice that the determinant of the matrix above is given by
$$n^2\sin^2 n(\pi-\alpha)+ n^2(\cos\alpha+\cos n(\pi-\alpha))^2.$$
This determinant is never zero, for if $\sin n(\pi-\alpha)=0$,
then $n(\pi-\alpha)=k\pi$, for some integer $k$. Nonetheless, in
this case $\cos n(\pi-\alpha)=(-1)^k$ and since $\alpha\in
(0,\pi)$, it is impossible to have $\cos\alpha+\cos
n(\pi-\alpha)=0$. It follows that the only solutions for the
system of equations is $a_n=b_n=0$ for all $n$. Thus, the
solutions of the differential equation (\ref{ecdif}) are constant
functions and $K$ must be a disc centred at $O$. \fin

\begin{theorem}\label{cuerda_l(t)}
Let $K$ be a strictly convex body in the plane and let
$\alpha\in(0,\pi)$ be a fixed angle. Suppose $\lambda(t)=\lambda_0$, for every $t\in [0,2\pi]$, for a positive number $\lambda_0,$ and $K$ has rotational symmetry of angle $\pi-\alpha$ or $2\alpha$. Then $K_{\alpha}$ is a
circle. 
\end{theorem}

\emph{Proof.} The vector $\nu(t)=\gamma_{\alpha}(t+\pi-\alpha)-\gamma_{\alpha}(t-\pi+\alpha)$ can be expressed using the parametrization given in (\ref{eq2}) by
\begin{align*}
\nu(t) &=[2p(t+\pi-\alpha)\cos \alpha +p(t-2\alpha)+p(t)]u(t)\\
&+\frac{1}{\sin\alpha}[p(t+\pi-\alpha)\cos 2\alpha -p(t-\pi+\alpha)+p(t-2\alpha)\cos \alpha-p(t)\cos\alpha]u'(t).
\end{align*}
If $K$ has rotational symmetry of angle $\pi-\alpha$ or $2\alpha$ we have that 
\begin{align*}
\nu(t) &=[2p(t+\pi-\alpha)\cos \alpha +2p(t)]u(t)\\
&+\frac{1}{\sin\alpha}[p(t+\pi-\alpha)\cos 2\alpha -p(t-\pi+\alpha)]u'(t),
\end{align*}
or equivalently $$\nu(t) =[2p(t+\pi-\alpha)\cos \alpha +2p(t)]u(t)+[p(t+\pi-\alpha)(-2\sin\alpha)]u'(t).$$ In this case $|\nu(t)|^2=\lambda^2(t)$ can be easily calculated as $$\lambda_0^2=4[p^2(t+\pi-\alpha)+p^2(t)+2p(t+\pi-\alpha)p(t)\cos\alpha],$$ hence $$\frac{\lambda_0^2}{4}=p^2(t+\pi-\alpha)+p^2(t)-2p(t+\pi-\alpha)p(t)\cos(\pi-\alpha)=d(t)^2,$$ i.e., the value of $d(t)=\frac{\lambda_0}{2}.$

\begin{figure}[H]
    \centering
    \includegraphics[width=.85\textwidth]{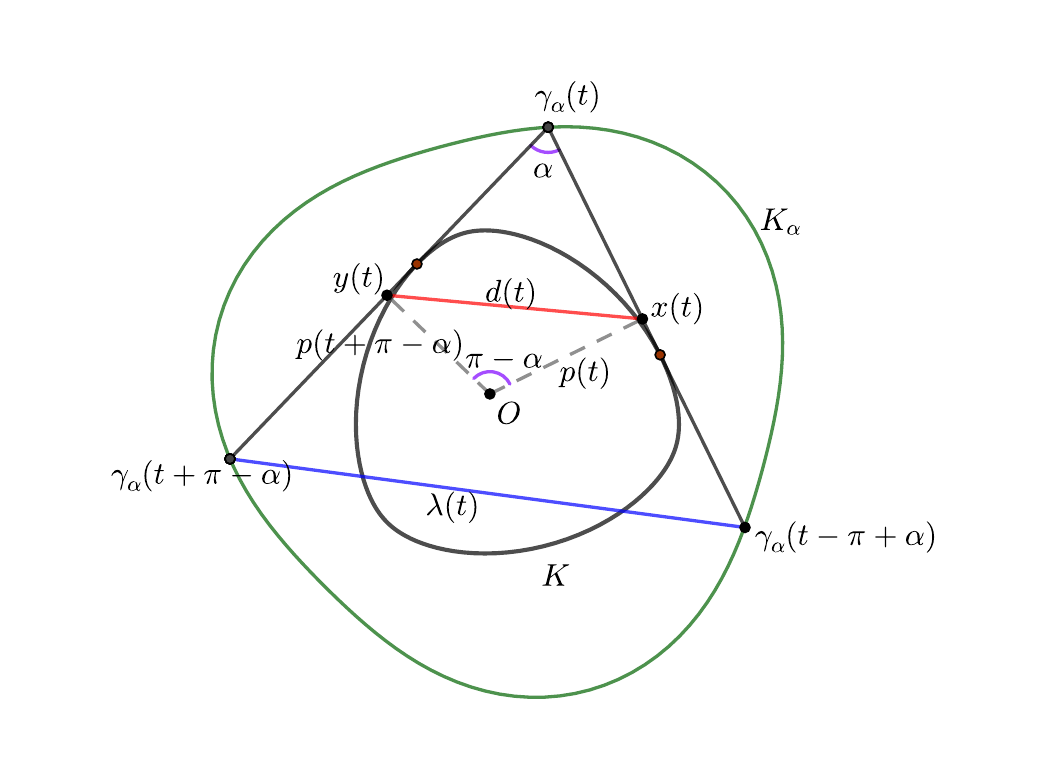}
    \caption{The value of $d(t)$ is constant}
    \label{parametro_d}
\end{figure}

Let $x(t)$ and $y(t)$ be the projections of $O$ into the support lines of $K$ at $\gamma(t)$ and $\gamma(t+\pi-\alpha)$, respectively (see Fig. \ref{parametro_d}). The quadrilateral $Ox(t)\gamma_{\alpha}(t)y(t)$ is cyclic, i.e., there exist a circle which passes through its four vertices, hence $|\gamma_{\alpha}(t)-O|=\frac{d(t)}{\sin\alpha}=\frac{\lambda_0}{2\sin\alpha}.$ It follows that $K_{\alpha}$ is a circle centred at $O$. \fin

\section{Some comments about rotors in the plane}
In this section we give some words about how the results obtained in this work are related to rotors in polygons. Moreover, in all the examples shown below, if we fix the convex body and the circumscribed polygon is rotated, while maintained circumscribed to $K$, the vertices describe a isoptic of $K$.

\begin{figure}[H]
    \centering
    \includegraphics[width=.41\textwidth]{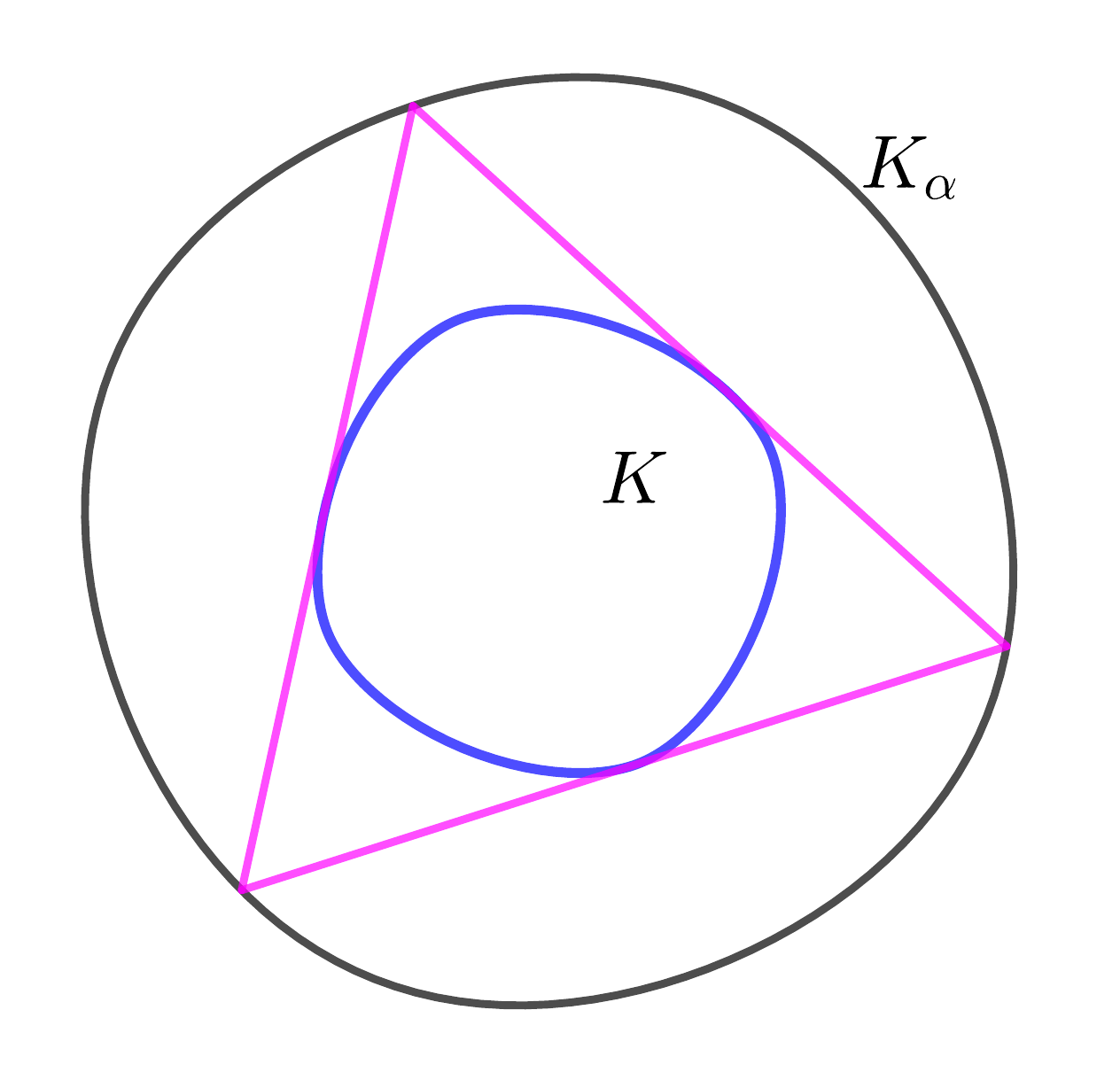}
    \caption{A centrally symmetric rotor in the equilateral triangle}
    \label{n_4}
\end{figure}

When $c(t)$ has constant value, using the Fourier series for the support function of $p$ in the proof of theorem \ref{cuerdas} we arrived to the equations
\begin{equation*}
-2\sin n(\pi-\alpha)\cos\alpha +\sin 2n\alpha=0
\end{equation*}
and 
\begin{equation*}
2\cos n(\pi-\alpha)\cos\alpha+1 +\cos 2n\alpha=0.
\end{equation*}
If $n$ is even, then
\begin{equation*}
\sin n\alpha (\cos n\alpha +\cos\alpha)=0
\end{equation*}
and 
\begin{equation*}
\cos n\alpha(\cos n\alpha+\cos\alpha)=0.
\end{equation*}
Both of them are zero if $\cos n\alpha+\cos\alpha=0,$ or after some trigonometric transformations $$\cos \left ( \frac{n\alpha +\alpha}{2}\right)\cos \left ( \frac{n\alpha -\alpha}{2}\right)=0.$$
It follows that $\alpha=\left( \frac{2r+1}{n+1}\right)\pi$ or $\alpha=\left( \frac{2r+1}{n-1}\right)\pi$, where $r$ is any integer number. In this case the determinant of the associated matrix is zero and we can choose the coefficients $a_n$, $b_n$, arbitrarily. For instance, for $n=4$ we select $\alpha=\frac{\pi}{3}$, $a_0=30$, $a_4=0$, $b_4=1$, and for any other natural number $n$ we have that $a_n=b_n=0$, i.e., the support function of $K$ is $p(t)=30+\sin 4t$ (see Fig. \ref{n_4}). The body $K$ in this example is centrally symmetric.

\begin{figure}[H]
    \centering
    \includegraphics[width=.46\textwidth]{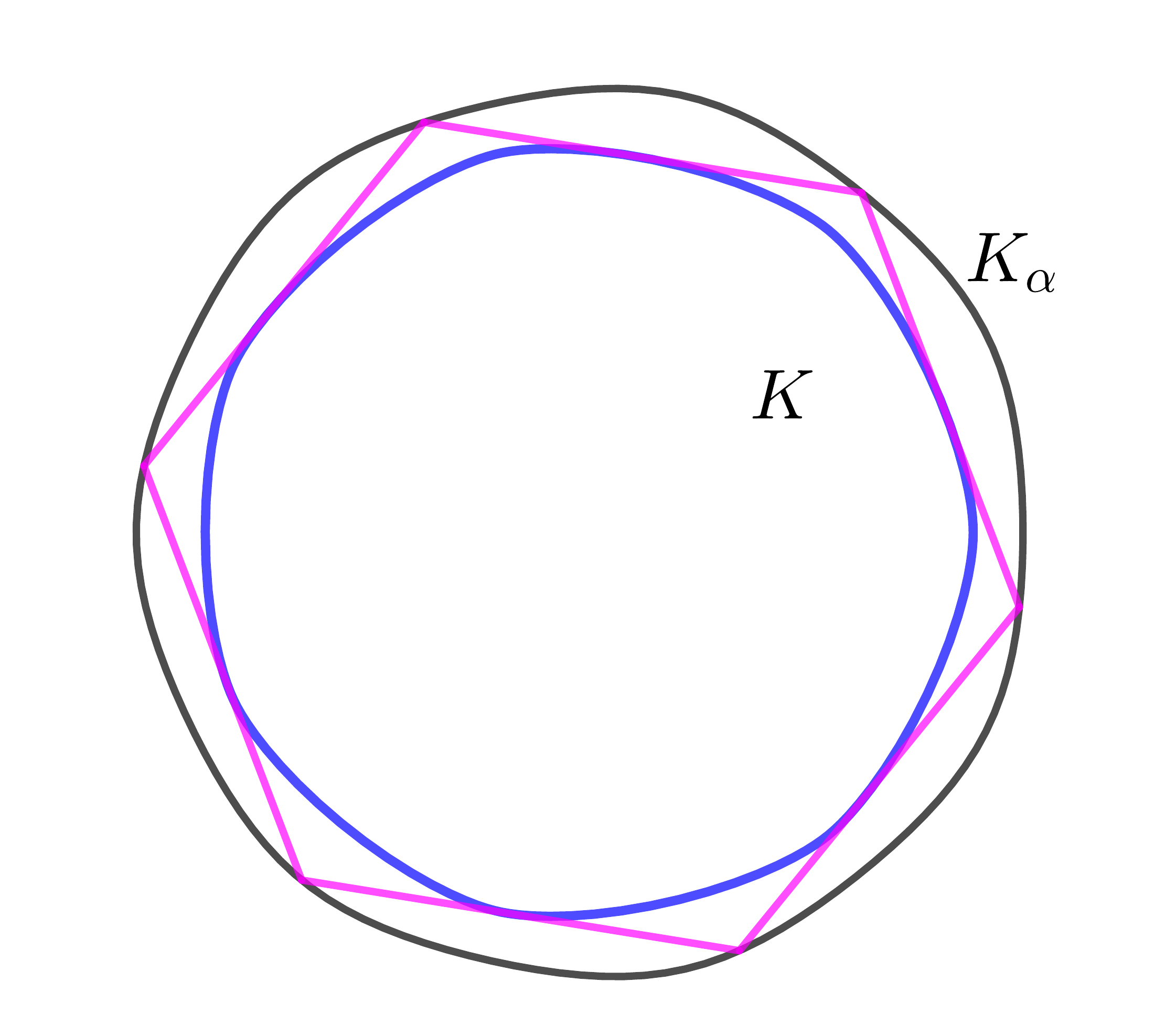}
    \caption{A rotor with constant width in the regular hexagon}
    \label{n11}
\end{figure}

If $n$ is odd, then
\begin{equation*}
\sin n\alpha (\cos n\alpha -\cos\alpha)=0
\end{equation*}
and 
\begin{equation*}
\cos n\alpha(\cos n\alpha-\cos\alpha)=0.
\end{equation*}
Both of them are zero if $\cos n\alpha-\cos\alpha=0,$ or after some trigonometric transformations $$\sin \left ( \frac{n\alpha +\alpha}{2}\right)\sin \left ( \frac{n\alpha -\alpha}{2}\right)=0.$$ 
It follows that $\alpha=\left( \frac{2r}{n+1}\right)\pi$ or $\alpha=\left( \frac{2r}{n-1}\right)\pi$, where $r$ is any integer number. Notice that we can obtain an example for any angle of the form $\alpha=\frac{s}{q}\pi$, where $s$ and $q$ are integers such that $0<s<q$. We just use this case, since $\alpha=\frac{2s}{2q}\pi$, we select $r=s$ and $n=2q+1$ or $n=2q-1$. For instance, for $\alpha=\frac{2\pi}{3}\pi=\frac{4\pi}{6}\pi$, and then we select $n=7$. In Fig. \ref{n11} we show the body $K$ with its isoptic $K_{2\pi/3},$ with the property that $c(t)$ has a constant value. The support function for this example is $p(t)=80+\cos 7t$. The body in this example has constant width.
 
\begin{figure}[H]
    \centering
    \includegraphics[width=.8\textwidth]{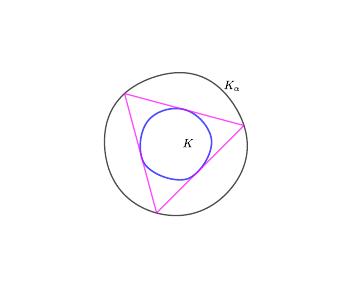}
    \caption{A rotor in the equilateral triangle, which is neither of constant width nor centrally symmetric}
    \label{n_4_5}
\end{figure}

The example shown in Fig. \ref{n_4_5} has coefficients different from zero for $n=4$ and $n=5$, consequently, the body $K$ obtained is neither of constant width or centrally symmetric. The support function of $K$ is $p(t)=70+\sin 4t+\cos 5t.$
 
Finally, we give an example of a rotor in the square. The support function is $p(t)=60+\cos 5t$.

\begin{figure}[H]
    \centering
    \includegraphics[width=.38\textwidth]{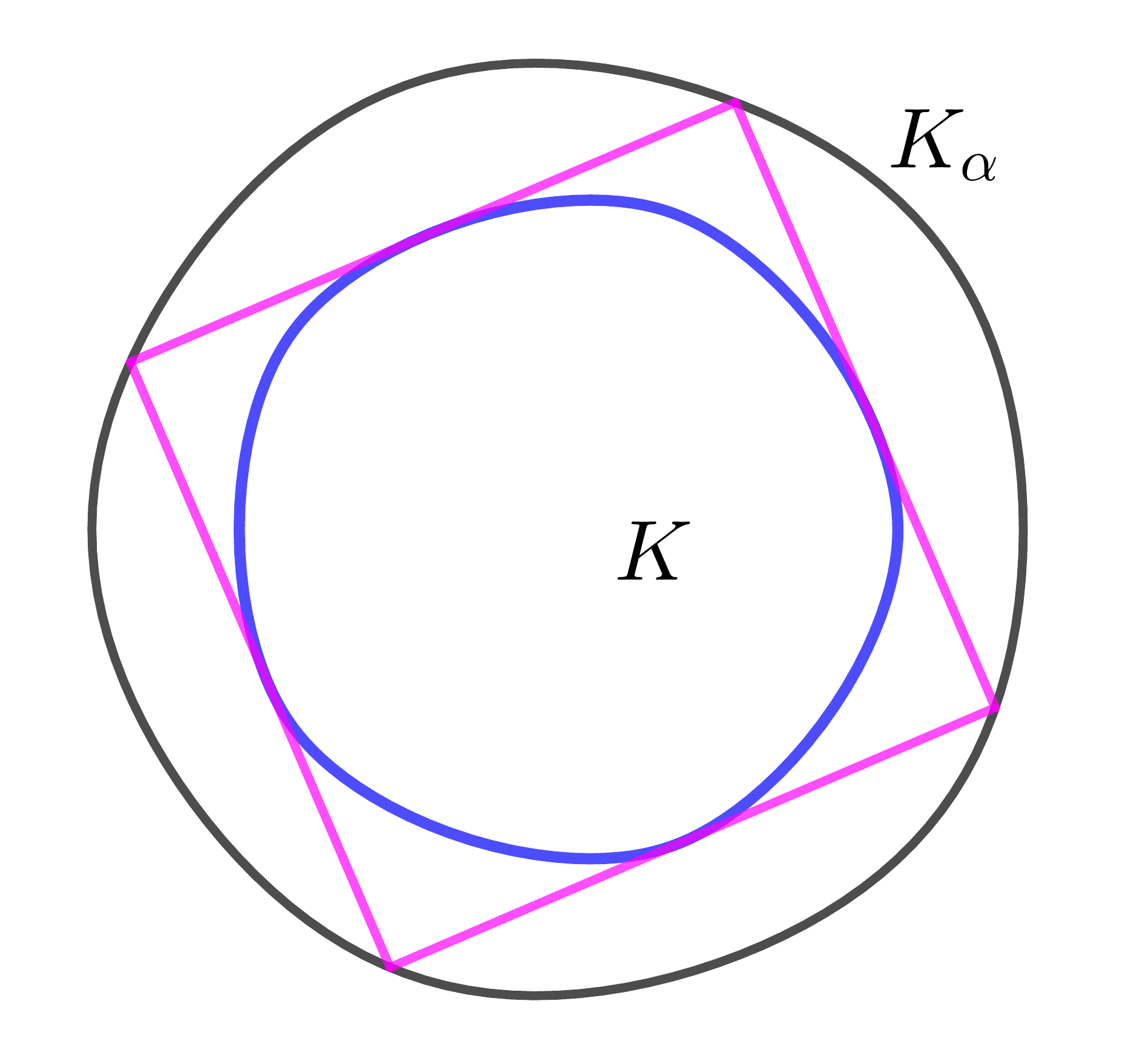}
    \caption{A rotor for the square}
    \label{pi_2}
\end{figure}

\section{An inequality about the length of some special chords}
The following inequality and characterization of the disc, in terms of the length $q(t)$, was proved in \cite{JeronimoYee}.

\textbf{Theorem JY.} Let $K$ be a strictly convex body in the plane with minimal width $\omega_0$. For any fixed $\alpha\in(0,\pi)$ there exists at least one value of the parameter $t$, $t(\alpha)\in [0,2\pi]$, such that $q(t(\alpha))\geq \omega_0\cos\frac{\alpha}{2}$. Moreover, if there is not such a chord with length exceeding $\omega_0\cos\frac{\alpha}{2}$, then $K$ is a disc.

Here we prove the following similar result.

\begin{theorem}\label{desigualdad}
Let $K$ be a convex body in the plane with minimal width $\omega_0$. For any fixed $\alpha\in(0,\pi)$ there exists at least one value of the parameter  $t$, $t(\alpha)\in [0,2\pi]$, such that $\lambda(t(\alpha))\geq 2\omega_0\cos\frac{\alpha}{2}$. Moreover, if there is not such a chord with length exceeding $2\omega_0\cos\frac{\alpha}{2}$, then $K$ is a convex body of constant width $\omega_0$.
\end{theorem}

\begin{figure}[H]
    \centering
    \includegraphics[width=.62\textwidth]{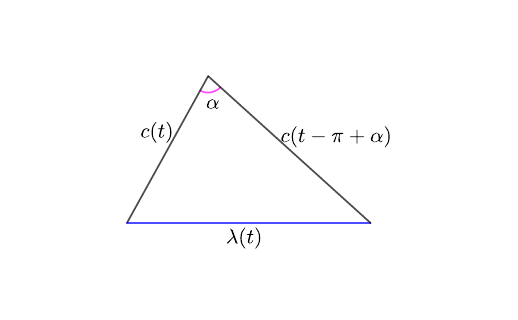}
    \caption{The minimum of $\lambda(t)$ is when $c(t)=c(t-\pi+\alpha)$}
    \label{triangulo}
\end{figure}

\emph{Proof.} The mean value of the chords of $K_{\alpha}$ that are tangent to $K$ is 
\begin{align*}
\overline{c(t)}&=\frac{1}{2\pi\sin\alpha}\int_0^{2\pi}[2p(t+\pi-\alpha)\cos\alpha+p(t)+p(t-2\alpha)]dt,\\
&=\frac{1}{2\pi\sin\alpha}[2L(K)\cos\alpha +2L(K)], \hspace{.6 cm}\text{(using Cauchy's formula (\ref{Cauchy}))}\\
&=\frac{\cot\frac{\alpha}{2}}{\pi}L(K).
\end{align*}
Indeed, by the same argument, the mean value of $c(t)+c(t-\pi+\alpha)$ is
$$\overline{c(t)+c(t-\pi+\alpha)}=\frac{2\cot\frac{\alpha}{2}}{\pi}L(K).$$
Let $t\in[0,2\pi]$ be any number and consider the triangle with sides of length $c(t)$, $c(t-\pi+\alpha)$, $\lambda(t)$, and angle $\alpha$, as shown in Fig. \ref{triangulo}. By Problem C4 in the book by I. Niven \cite{Niven}, we have that the minimum of $\lambda(t)$, among all triangles with a fixed angle $\alpha$ and the sum $c(t)+ c(t-\pi+\alpha)=c_0$, for a constant number $c_0$, is obtained if and only if $c(t)=c(t-\pi+\alpha)=c_0/2$. It follows that for every $t\in[0,2\pi]$ it holds that
$$\lambda(t)\geq [c(t)+c(t-\pi+\alpha)]\sin\frac{\alpha}{2},$$ with equality if and only if $c(t)=c(t-\pi+\alpha).$ Let $t_0\in[0,2\pi]$ be such that $c(t_0)+c(t_0-\pi+\alpha)=\overline{c(t)+c(t-\pi+\alpha)},$ we have that 
\begin{equation}\label{desi}
\lambda(t_0)\geq [c(t_0)+c(t_0-\pi+\alpha)]\sin\frac{\alpha}{2}=\frac{2\cos\frac{\alpha}{2}}{\pi}L(K).
\end{equation} 

Using Cauchy's formula for the perimeter of $K$ we can also prove that $L(K)\geq \pi \omega_0$, with equality if and only $K$ is a body of constant width $\omega_0.$ Hence $\lambda(t_0)\geq 2\omega_0\cos\frac{\alpha}{2}.$

\begin{figure}[H]
    \centering
    \includegraphics[width=.45\textwidth]{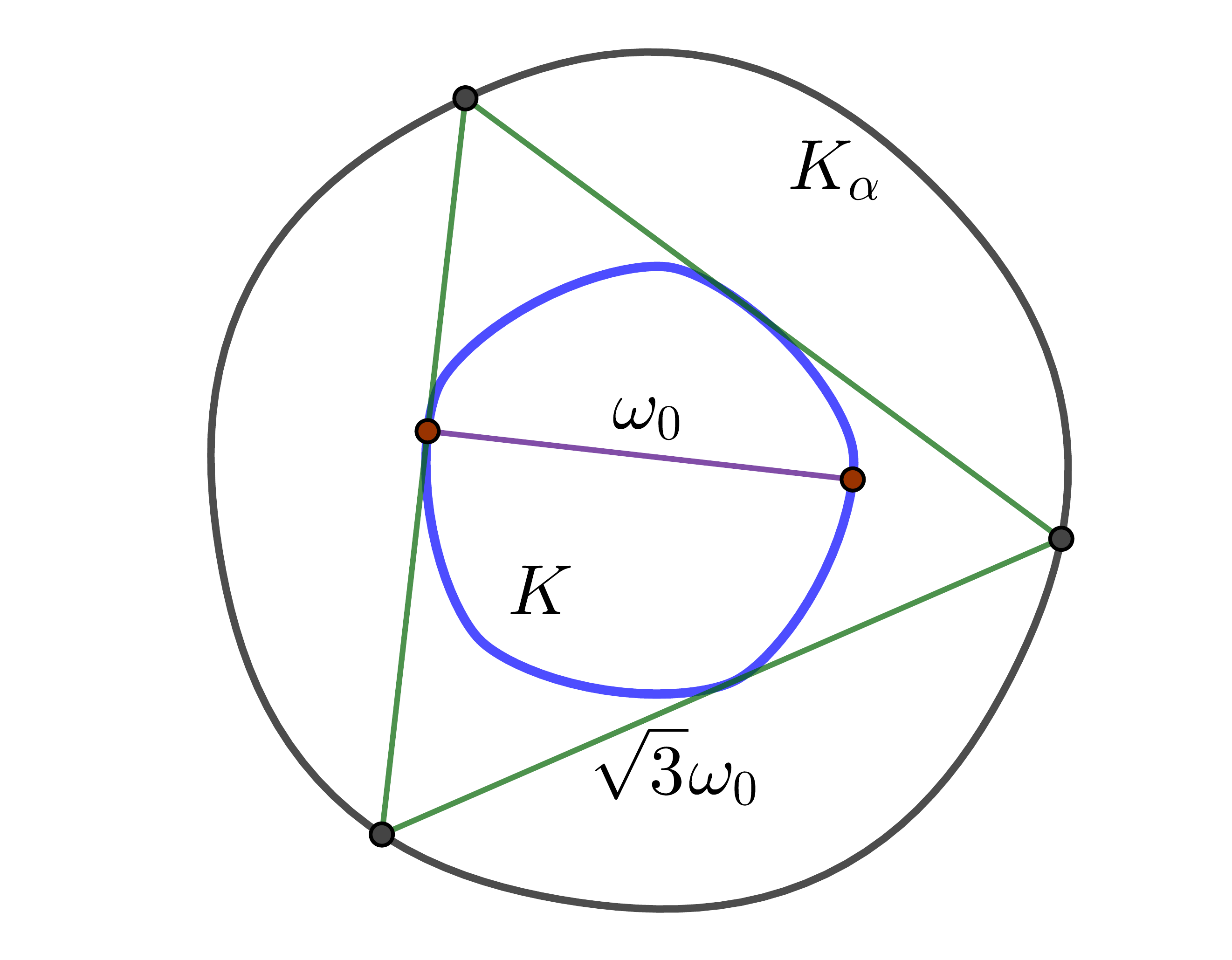}
    \caption{A convex body with $\lambda(t)=2\omega_0\cos\frac{\alpha}{2}$ for every $t$}
    \label{cte}
\end{figure}

Now, if there is no chord with $\lambda(t)>2\omega_0\cos\frac{\alpha}{2}$, then $\lambda(t_0)=2\omega_0\cos\frac{\alpha}{2}$, which by (\ref{desi}) implies that $L(K)=\pi \omega_0,$ i.e., $K$ must be a body of constant width $\omega_0$ (see Fig. \ref{cte}). \fin

%%%%%%%%%%%%%%%%%%%%%%%%%%%%%%%%%%%%%%%%%%%%%%%%%%%%%%%%%%%%%%%%%%%%%%%%%%%%%%%%%%
\section{Some results in higher dimensions}
The following lemma is needed for some of the subsequent results.

\begin{lemma}\label{circulos}
Let $K,L\subset \mathbb R^n$, $n\geq 3$, be convex bodies with $L\subset \emph{int}\, K$, such that every $(n-1)$-dimensional section of $K$ tangent to $L$ is an $(n-1)$-dimensional ball, then $K$ is a ball. If additionally, the centre of every tangent ball is at the point of contact with $L$, then $K$ and $L$ are concentric balls.
\end{lemma}

\emph{Proof.} There are several ways to prove this lemma. Here we give the following proof. First note that if all the $3$-dimensional sections of a convex body through a given point in its interior are $3$-dimensional balls, then it is an $n$-dimensional ball. We consider any point in the interior of $L$ for such a point and since the hypothesis of the theorem is inheritated to every $3$-dimensional section, we have that it is sufficient to prove the theorem in the $3$-dimensional case.

Let $x$ be any point in $\partial K$ and let $\ell$ be any line supporting $K$ at $x$. Consider the two support planes of $L$ which share the line $\ell$, to say $H_1$ and $H_2$. Let $H$ be any other support plane of $L$ through the point $x$. By hypothesis $H\cap K$ is a disc which intersects each one of the circles $H_1\cap\partial K$ and $H_2\cap \partial K$ at two points. The circle $H\cap \partial K$ passes through three points of the set $(H_1\cap \partial K)\cup (H_2\cap \partial K)$, and since there is a unique sphere which contains the two circles $H_1\cap \partial K$  and $H_2\cap \partial K$, it holds that this sphere contains $H\cap\partial K$. Since $H$ is any support plane of $L$ through $x$, we have that $\mathcal R_x=\{H\cap\partial K: H \text{ is a support plane of } L \text{ through } x\}$ is a closed subset of a sphere. Let $y\in \partial K$ be any point such that $\mathcal R_y\cap \mathcal R_x\neq\emptyset$, then $\mathcal R_y\cup \mathcal R_x$ is contained in the same sphere. Continuing in this way, since $K$ is a compact set, we can prove that $\partial K$ is a sphere.

Now, let $H$ be any hyperplane supporting $L$ at a point $z$ and suppose the centre of the $(n-1)$-dimensional ball $H\cap K$ is $z$. The line orthogonal to $H$ through $z$ passes through the centre of the ball $K$, this implies that indeed $H$ is the tangent plane of $L$ through $z$, i.e., $\partial K$  is a differentiable surface. We have that all the normal lines of $\partial L$ passes through the centre of $K$, this implies (see for instance \cite{Toponogov}) that $\partial L$ is a sphere with centre at the centre of $K$. We conclude that $K$ and $L$ are concentric balls. \fin

Let $\beta=\pi-\alpha$. In \cite{Green},  J. W.
Green proved that $K$ is a Euclidean disc if $K_{\alpha}$ is a
circle and any of the following conditions hold: $\beta$ is an
irrational multiple of $\pi$, or $\beta=\frac{2m}{n}\pi$, with
$2m$ and $n$ relatively prime integer numbers. On the other side,
M.S. Klamkin conjectured \cite{Klamkin} and J. Nitsche proved
\cite{Nitsche} that two different values $\alpha_1$, $\alpha_2$,
such that $K_{\alpha_1}$ and $K_{\alpha_2}$ are circles, are
enough to prove that $K$ is a Euclidean disc. A similar result was proved by \'A. Kurusa and T. \'Odor in \cite{Kurusa_Odor} in $\mathbb R^n$, $n\geq 3$, where the isoptic surface is defined using the solid angle under which is seen a convex body $K$. They proved that if two isoptic surfaces of a given convex body are concentric spheres, then it is a ball. The following is a variant of this kind of results.

\begin{theorem}\label{secciones}
Let $K\subset \mathbb R^3$ be a strictly convex body contained in the interior of a sphere $\mathcal S$. Suppose that for every $2$-dimensional plane $H$, with $H\cap K\neq \emptyset$, it holds that $H\cap\mathcal S$ is a isoptic curve of $H\cap K.$ Then $K$ is a ball concentric with $\mathcal S$.
\end{theorem}

\emph{Proof.} We first observe the following: \emph{If for some angle $\alpha\in (0,\pi)$ the $\alpha$-isoptic curve of a planar body $M$ is a circle $S$, then its centre $x$ belongs to $M$.} Suppose this is not the case. Let $y\in M$ be the point which is closest to $x$. Let $z\in \mathcal S$ be the point where the ray $\overrightarrow{xy}$ intersects to $\mathcal S$ and let $z'$ be the point in $\mathcal S$ such that $[z,z']$ is a diameter. Let $\ell$ be the line through $x$ which is orthogonal to $[z,z']$. The central projection of $M$ into the line $\ell$  is a segment $[a,b]$ that contains the centre $x$. The angle formed by the lines $z'a$ and $z'b$ contains $M$, however, at least one of the lines $z'a$ or $z'b$ does not intersect $M$. Since the angles $\measuredangle azb$ and $\measuredangle az'b$ are equal, we obtain that the angles under which $M$ is seen from the points $z$ and $z'$ are not equal. This contradicts that $S$ is the isoptic curve of angle $\alpha$ of $M$. Hence, we have that $x$ belongs to $M$.

\begin{figure}[H]
    \centering
    \includegraphics[width=.5\textwidth]{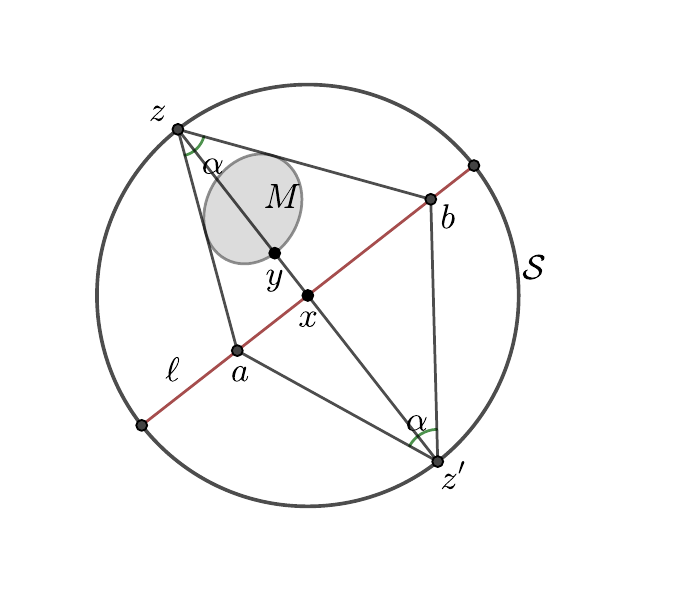}
    \caption{The centre $x$ belong to $M$}
    \label{circulo}
\end{figure}

Now let $q\in\partial K$ be any arbitrary point and let $H_q$ be a plane supporting $K$ at $q$. By hypothesis we have that $H_q\cap \mathcal S$ is a circle; we are going to prove that its centre is $q$. Consider a sequence of planes $\{H_r\}$ parallel to $H_q$ in such a way that the sequence converges to $H_q$. Let $x_r$ be the centre of the circle $H_r\cap \mathcal S$, for every natural number $r$. By the comment above, we have that $x_r\in H_r\cap K$, for every $r\in\mathbb N$, and since $H_r\cap K\longrightarrow H_q\cap K=q$ we have that $x_r\longrightarrow q$ when $r\longrightarrow \infty$. In other words, the centre of $H_q\cap \mathcal S$ is $q$. Now we apply Lemma \ref{circulos} and conclude that $K$ is a ball concentric with $\mathcal S$. \fin

The following result shows that in dimension 3 or higher, only Euclidean balls have an equichordal convex body.

\begin{theorem}\label{equichordal}
Let $K\subset \mathbb R^n$, $n\geq 3$, be a strictly convex body which possesses an equichordal convex body $L$ in its interior. Then $K$ and $L$ are concentric balls.
\end{theorem}

\emph{Proof.} Suppose the length of the chords of $K$ tangent to $L$ is $\lambda$. First we prove that $L$ is strictly convex, i.e., $L$ has not segments in its boundary. Suppose to the contrary that there is a segment $[a,b]\subset \partial L$ and consider any $2$-dimensional plane $H$ supporting $L$ at $[a,b]$. By hypothesis, every chord of $H\cap K$ through $a$ has length $\lambda$, and the same happens for every chord through $b$. This implies that $H\cap K$ has two equichordal points, but as were proved by M. R. Rychlik \cite{Rychlik} this is not possible.

Now, let $H$ be any $2$-dimensional plane supporting $L$ at a point $x$. We will prove that the diameter of $H\cap K$ is $\lambda$. Suppose this is not the case and there is a chord $[a,b]$ of $H\cap K$ with $|a-b|>\lambda.$ Clearly, $[a,b]\cap L=\emptyset.$ Consider a $2$-dimensional plane $\Pi$ that contains to $[a,b]$ and such that $\Pi\cap \text{int}\, L \neq\emptyset.$ Let $c,d\in\partial (\Pi\cap L)$ such that the lines through $c$ and $d$, respectively, that are parallel to $[a,b]$ are support lines of $\Pi\cap L$. Since $|a-b|>\lambda$, we have that there exist a chord $[e,f]$ of $\Pi\cap K$ that is separated from $\Pi\cap L$ by the chord $[a,b]$ and $|e-f|=\lambda$. We have three chords of $\Pi\cap K$ parallel to $[a,b]$ and with length $\lambda$. This contradicts that $K$ is strictly convex, then the diameter of $H\cap K$ must be $\lambda$. 

Every chord of $H\cap K$ through $x$ has length equal to $\lambda$, then every chord of $H\cap K$ through $x$ is a binormal of $H\cap K$; it follows that $H\cap K$ is a disc centred at $x$ (see for instance \cite{Chakerian_Groemer}). Since this is true for every $2$-dimensional plane $H$ through $x$, we have that any $(n-1)$-dimensional section of $K$ tangent to $L$ at $x$ is an $(n-1)$-dimensional sphere centred at $x$. This is also true for every $x\in\partial L$, so we apply Lemma \ref{circulos} and conclude that $K$ and $L$ are concentric balls. \fin

\begin{corollary}\label{tajadas}
Let $K,L\subset \mathbb R^n$, $n\geq 3$, be two convex bodies. Suppose that for every hyperplane $H$ that intersects to $L$ it holds that $H\cap L$ is an equichordal body of $H\cap K$. Then $K$ and $L$ are concentric balls.
\end{corollary}

\emph{Proof.} It is quite simple to prove that all the chords of $K$ tangent to $L$ have the same length. The conclusion follows applying Theorem \ref{equichordal}. \fin

Suppose now that $K$ is a strictly convex body in $\mathbb R^3$. For $z_1,z_2\in\partial K$, we will say that the segment $[z_1,z_2]$ is an $\alpha$-chord of $K$, if there exist tangent planes at $z_1$ and $z_2$ meeting at an angle $\alpha$. We present the following analogue to Theorem JY.

\begin{theorem}\label{cuerdas_dim3}
Let $K$ be a strictly convex body in $\mathbb R^3$ with boundary of class $C^1$ and with minimal width $\omega_0$. For any fixed $\alpha\in(0,\pi)$ there exists an $\alpha$-chord of $K$ with length at least $\omega_0\cos\frac{\alpha}{2}$. Moreover, if there is not an $\alpha$-chord with length exceeding $\omega_0\cos\frac{\alpha}{2}$, then $K$ is a ball.
\end{theorem}

\emph{Proof.} Let $u_0\in\mathbb S^2$ be  a direction such that the width of $K$ in this direction has minimal value $\omega_0$. Let $\upsilon\in\mathbb S^2$ be a vector orthogonal to $u_0$ and let's denote the orthogonal
projection on the plane $\upsilon^{\bot}$ by $\pi_{\upsilon}$. Clearly, the minimum width of $K'=\pi_{\upsilon}(K)$ is also
$\omega_0$ and every $\alpha$-chord of $K'$ is obtained as the projection of some $\alpha$-chord of $K$. By Theorem JY there is an $\alpha$-chord of $K'$ with length at least $\omega_0\cos\frac{\alpha}{2}$ and so there is an
$\alpha$-chord of $K$ with length larger than or equal to $\omega_0\cos\frac{\alpha}{2}$.

Suppose now that all the $\alpha$-chords of $K$ have length less than or equal to $\omega_0\cos\frac{\alpha}{2}$. Then all the
$\alpha$-chords of $\pi_{\upsilon}(K)$, for every $\upsilon\in\mathbb S^2$, have length less than or equal to $\omega_0\cos\frac{\alpha}{2}$. It follows from Theorem JY that $\pi_{\upsilon}(K)$ is a disc. Since $\upsilon$ is an arbitrary direction, we have that every $2$-dimensional projection of $K$ is a disc and therefore $K$ is a ball. \fin

\begin{theorem}\label{cuerdas_cte_dim3}
Let $K$ be a strictly convex body in $\mathbb R^3$ with boundary of class $C^1$ and with $\frac{\pi}{2}$-chords of constant length $\lambda$, then $K$ is a ball.
\end{theorem}

\emph{Proof.} For any point $z\in\partial K$, let $H_z$ denote the tangent plane of $K$ at $z$ and let $A_z=\{w\in\partial K:[z,w]\text{ is an }\frac{\pi}{2}\text{-chord of }K\}$. For each point $w$ in $A_z$, we have that $H_w$ intersects
$H_z$ with angle $\frac{\pi}{2}$. Let $z,z'\in\partial K$ be two points such that $[z,z']$ is normal to $H_z$ and $H_{z'}$. For every $w\in A_z$, the triangle $\Delta zz'w$ is isosceles with sides $[z,w]$ and $[z',w]$ of equal length. Then $A_z$ is a circle centred at $(z+z')/2$. Moreover, given any two antipodal points $z_1$ and $z_1'$ of $A_z$, the segment $[z_1,z_1']$ must be orthogonal to the planes $H_{z_1}$ and $H_{z_1'}$. By  repeating the previous argument to each pair of antipodal points of $A_z$, we have that any section of $K$ containing $[z,z']$ is a disc. This implies that $K$ is a ball. \fin

This manuscript has not any associated data.

\end{document}